\documentclass[12]{amsart}
\usepackage{amsmath,amssymb,amsthm,color,enumerate,comment,centernot,enumitem,url,cite}
\usepackage{graphicx,relsize,bm}
\usepackage{mathtools}
\usepackage{array}

\makeatletter
\newcommand{\tpmod}[1]{{\@displayfalse\pmod{#1}}}
\makeatother

\newtheorem{thm}{Theorem}[section]

\newtheorem{cor}[thm]{Corollary}

\theoremstyle{remark}

\theoremstyle{definition}

\newtheorem{rem}[thm]{Remark}

\theoremstyle{THM}

\newcommand{\abs}[1]{\left|{#1}\right|}

\def\FF {{\mathcal F}}
\def\GG {{\mathcal G}}

\def\Z {{\mathbb Z}}

\def\Q {{\mathbb Q}}

\def\D {{\mathcal D}}

\def\GG {{\mathcal G}}

\def\D {{\mathcal D}}
\def\Z {{\mathbb Z}}
\def\Q {{\mathbb Q}}

\def\Gal{{\mbox{{\rm{Gal}}}}}

\makeatletter
\@namedef{subjclassname@2020}{%
  \textup{2020} Mathematics Subject Classification}
\makeatother

\def\red#1 {\textcolor{red}{#1 }}
\def\blue#1 {\textcolor{blue}{#1 }}

\numberwithin{equation}{section}

\def\Z {{\mathbb Z}}

\newcommand{\Mod}[1]{\ (\mathrm{mod}\enspace #1)}

\begin{document}

\title[Monogenic $x^4+ax^3+d$ and their Galois groups]{Monogenic trinomials of the form $x^4+ax^3+d$ and their Galois groups}

\author{Joshua Harrington}
\address{Department of Mathematics, Cedar Crest College, Allentown, Pennsylvania, USA}
\email[Joshua Harrington]{Joshua.Harrington@cedarcrest.edu}

\author{Lenny Jones}
\address{Professor Emeritus, Department of Mathematics, Shippensburg University, Shippensburg, Pennsylvania 17257, USA}
\email[Lenny~Jones]{doctorlennyjones@gmail.com}

\date{\today}

\begin{abstract}
Let $f(x)=x^4+ax^3+d\in \Z[x]$, where $ad\ne 0$. Let $C_n$ denote the cyclic group of order $n$, $D_4$ the dihedral group of order 8, and $A_4$ the alternating group of order 12. Assuming that $f(x)$ is monogenic, we give necessary and sufficient conditions involving only $a$ and $d$ to determine the Galois group $G$ of $f(x)$ over $\Q$. In particular, we show that $G=D_4$ if and only if $(a,d)=(\pm 2,2)$, and that $G\not \in \{C_4,C_2\times C_2\}$. Furthermore, we prove that $f(x)$ is monogenic with $G=A_4$ if and only if $a=4k$ and $d=27k^4+1$, where $k\ne 0$ is an integer such that $27k^4+1$ is squarefree.  This article extends previous work of the authors on the monogenicity of quartic polynomials and their Galois groups.
\end{abstract}


\maketitle
\section{Introduction}\label{Section:Intro}

  Throughout this article, we let $a,d\in \Z$ with $ad\ne 0$, and we let
  \begin{equation}\label{Eq:frdelta}
  f(x):=x^4+ax^3+d, \quad r(x):=x^3-4dx-a^2d\quad \mbox{and} \quad \delta:=256d-27a^4.
  \end{equation}
    We let $\Gal(f)$ denote the Galois group of $f(x)$ over $\Q$, when $f(x)$ is irreducible over $\Q$.  We use Maple to calculate
 \begin{equation}\label{Eq:Delta(f)}
 \Delta(f)=\Delta(r)=d^2\left(256d-27a^4\right)=d^2\delta,
 \end{equation}
 where $\Delta(f)$ and $\Delta(r)$ are the polynomial discriminants over $\Q$ of $f(x)$ and $r(x)$, respectively.

  We say that $f(x)$ is \emph{monogenic} if $f(x)$ is irreducible over $\Q$ and $\{1,\theta,\theta^2,\theta^3\}$ is a basis for the ring of integers $\Z_K$ of $K=\Q(\theta)$, where $f(\theta)=0$. Hence, $f(x)$ is monogenic if and only if  $\Z_K=\Z[\theta]$. For $f(x)$ irreducible over $\Q$, it is well known \cite{Cohen} that
\begin{equation} \label{Eq:Dis-Dis}
\Delta(f)=\left[\Z_K:\Z[\theta]\right]^2\Delta(K),
\end{equation}
where $\Delta(K)$ is the discriminant over $\Q$ of the number field $K$.
Thus, from \eqref{Eq:Dis-Dis}, $f(x)$ is monogenic if and only if $\Delta(f)=\Delta(K)$.

The possible Galois groups for a quartic polynomial are
\begin{equation}\label{Eq:GG}
C_4,\quad C_2\times C_2\quad D_4, \quad A_4, \quad S_4,
\end{equation}
where $C_n$ is the cyclic group of order $n$, $D_4$ is the dihedral group of order 8, $A_4$ is the alternating group of order 12, and $S_4$ is the symmetric group of order 24. With the exception of $C_4$, for each group $G$ in \eqref{Eq:GG}, the authors recently gave in \cite{HJquartic} an infinite family of monogenic quartic polynomials with Galois group $G$. In a private communication, Tristan Phillips subsequently asked the second author if it is possible to determine all monogenic quartic trinomials that have Galois group $C_4$. In \cite{JonesBAMSevenquartics}, the second author gave a partial answer to this question by showing that if $x^4+bx^2+d$ is monogenic with Galois group $C_4$, then $(b,d)\in \{(-4,2),(4,2),(-5,5)\}$.
The results in this article extend this previous work of the authors. More precisely, we prove the following.
\begin{thm}\label{Thm:Main1}
 If $f(x)$ is monogenic, then $\Gal(f)\simeq$
\begin{enumerate}
\item \label{Main1 I1} $D_4$ if and only if $(a,d)=(\pm 2,2)$; 
\item \label{Main1 I2} $A_4$ if and only if $\delta$ is a square;
\item \label{Main1 I3} $S_4$ if and only if $\delta$ is not a square and $(a,d)\ne (\pm 2,2)$.
\end{enumerate}
\end{thm}
 The following corollary is then immediate from Theorem \ref{Thm:Main1}.
 \begin{cor}\label{Cor:Main1}
  If $f(x)$ is monogenic, then
  \[\Gal(f)\not \simeq C_2\times C_2 \quad \mbox{and} \quad \Gal(f)\not \simeq C_4.\]
 \end{cor}
  \begin{thm}\label{Thm:Main2}
  The trinomial $f(x)$ is monogenic with $\Gal(f)\simeq A_4$ if and only if $a=4k$ and $d=27k^4+1$, where $k\ne 0$ is an integer such that  $d$ is squarefree.
\end{thm}

The following corollary then follows immediately from \cite{Granville}.
\begin{cor}\label{Cor:Main2}
  Let $k\ne 0$, and let $(a,d)=(4k,27k^4+1)$. Then, under the assumption of the $abc$-conjecture for number fields, there exist infinitely many values of $k$ such that $f(x)$ is monogenic and $\Gal(f)\simeq A_4$.
\end{cor}

We emphasize three important implications of Theorem \ref{Thm:Main1}, Corollary \ref{Cor:Main1} and Theorem \ref{Thm:Main2}. Firstly, we see that item \eqref{Main1 I1} of Theorem \ref{Thm:Main1} gives a result for $D_4$-trinomials $f(x)$ similar to the result for $C_4$-even trinomials found in \cite{JonesBAMSevenquartics}. Secondly, Corollary \ref{Cor:Main1} provides additional information toward a complete answer to the question of Phillips, since no trinomials $f(x)$ exist with $\Gal(f)\simeq C_4$. Thirdly, Theorem \ref{Thm:Main2} gives a complete and explicit
 description of all monogenic trinomials $f(x)$ with $\Gal(f)\simeq A_4$.

\section{Preliminaries}\label{Section:Prelim}
The first theorem below follows from a result due to Jakhar, Khanduja and Sangwan \cite[Theorem 1.1]{JKS2} for arbitrary irreducible trinomials when applied to our specific quartic trinomial $f(x)$.
\begin{thm}\label{Thm:JKS} 
Suppose that $f(x)$ is irreducible over $\Q$ with $f(\theta)=0$. Let $K=\Q(\theta)$, and let $\Z_K$ denote the ring of integers of $K$. Then $f(x)$ is monogenic if and only if every prime divisor $q$ of $\Delta(f)$ (in \eqref{Eq:Delta(f)}) satisfies one of the following conditions:
\begin{enumerate}[font=\normalfont]
  \item \label{JKS:I1} when $q\mid d$, then $q^2\nmid d$;
  \item \label{JKS:I2} when $q\mid a$ and $q\nmid d$, then
  \[\mbox{either} \quad q\mid a_2 \mbox{ and } q\nmid d_1 \quad \mbox{ or } \quad q\nmid a_2\left(-d^3a_2^4-d_1^4\right),\]
  where $a_2=a/q$ and $d_1=\frac{d+(-d)^{q^j}}{q}$, where $j$ is the highest power of $q$ such that $q^j$ divides $4$;
         \item \label{JKS:I3} when $q\nmid ad$, then $q^2\nmid \delta$.
   \end{enumerate}
\end{thm}

The following useful corollary is immediate from Theorem \ref{Thm:JKS}.
\begin{cor}\label{Cor:squarefree}
If $f(x)$ is monogenic, then $d$ is squarefree.
\end{cor}

The next theorem follows from a result due to Kappe and Warren \cite[Theorem 1]{KW} when applied to our specific quartic trinomial $f(x)$.
\begin{thm}\label{Thm:KW}
Suppose that $f(x)$ is irreducible over $\Q$. Let $L$ be the splitting field of $r(x)$, as defined in \eqref{Eq:frdelta}. Then $\Gal(f)\simeq$
\begin{enumerate}
 \item \label{GI:1} $C_4$ if and only if $r(x)$ has exactly one root $t\in \Z$ and
 \begin{equation}\label{Eq:gr}
 g(x):=(x^2-tx+d)(x^2+ax-t)
 \end{equation} splits over $L$;
 \item \label{GI:2} $C_2\times C_2$ if and only if $r(x)$ splits into linear factors over $\Z$;
 \item \label{GI:3} $D_4$ if and only if $r(x)$ has exactly one root $t\in \Z$ and $g(x)$, as defined in \eqref{Eq:gr}, does not split over $L$;
 \item \label{GI:4} $A_4$ if and only if $r(x)$ is irreducible over $\Z$ and $\delta$ is a square in $\Z$;
  \item \label{GI:5} $S_4$ if and only if $r(x)$ is irreducible over $\Z$ and $\delta$ is not a square in $\Z$.
 \end{enumerate}
 \end{thm}
 \begin{rem}
   The polynomial $r(x)$ in Theorem \ref{Thm:KW} is known as \emph{the cubic resolvent of $f(x)$}.
 \end{rem}

\section{The Proof of Theorem \ref{Thm:Main1}}\label{Section:Main1Proof}
\begin{proof}
Suppose first that $r(x)$ is reducible over $\Z$, so that $r(t)=0$ for some $t\in \Z$. Hence, for some $A,B\in \Z$, we have that
\begin{equation}\label{Eq:r2}
r(x)=(x-t)(x^2+Ax+B)=x^3+(A-t)x^2+(B-tA)x-tB.
\end{equation} Equating coefficients in \eqref{Eq:frdelta} and \eqref{Eq:r2} yields
\begin{equation}\label{Eq:r3}
 r(x)=(x-t)(x^2+tx+t^2-4d).
\end{equation} Then, calculating $\Delta(r)$  in  \eqref{Eq:r3}, and recalling $\Delta(r)$ from \eqref{Eq:Delta(f)}, yields the equation
\begin{equation}\label{Eq:Delta(r)}
\Delta(r)=d^2\delta=d^2(256d-27a^4)=(16d-3t^2)(4d-3t^2)^2.
\end{equation}
 Suppose that $d\not \in \D=\{1,2\}$. Since $r(t)=0$, we see from \eqref{Eq:frdelta} that $t^3=d(4t+a^2)$.  Therefore, since $d$ is squarefree from Corollary \ref{Cor:squarefree}, we deduce that $t=dw$ for some $w\in \Z$. Hence, since $\Delta(f)=\Delta(r)$, we have from \eqref{Eq:Delta(f)} and \eqref{Eq:Delta(r)} that
\begin{equation}\label{Eq:Del(f)}
\Delta(f)=d^2\delta=d^3(16-3dw^2)(4-3dw^2)^2.
\end{equation} 
Note that $4-3dw^2=-1$ is impossible, and $4-3dw^2=1$ implies that $d=1\in \D$. Hence, $\abs{4-3dw^2}>1$. Suppose that $q$ is an odd prime divisor of $4-3dw^2$. Then $q\nmid d$ and therefore, from \eqref{Eq:Del(f)}, we see that $q^2\mid \delta$. If $q\mid a$, then $q\mid 256d$, which is impossible. Hence, $q\nmid a$, 
which contradicts the fact that $f(x)$ is monogenic by item \eqref{JKS:I3} of Theorem \ref{Thm:JKS}. Thus, it follows that
\begin{equation}\label{Eq:powerof2}
4-3dw^2=\pm 2^m,
\end{equation} for some integer $m\ge 1$. Checking the cases $m\in \{1,2,3\}$, we arrive at
\[3dw^2\in\left\{\begin{array}{cl}
  \{2,6\} & \ \mbox{when $m=1$}\\[.5em]
  \{0,8\} & \ \mbox{when $m=2$,}\\[.5em]
  \{-4,12\} & \ \mbox{when $m=3$,}
\end{array}\right.\]
 which are impossible with the two exceptions $3dw^2\in \{6,12\}$. However, these exceptional cases imply, respectively, that $d=2$ and $d=1$, both of which are elements of $\D$. Thus, we have from \eqref{Eq:powerof2} that
\begin{equation}\label{Eq:dw}
3dw^2=4(1\pm 2^{m-2}),
\end{equation} for some integer $m\ge 4$. Note that \eqref{Eq:dw} implies that $2\nmid d$ since $d$ is squarefree, so that $w=2z$, for some odd integer $z$.  Recalling that $t=dw$ and $r(t)=0$, we have
\[r(t)=t^3-4dt-a^2d=d(8d^2z^3-8dz-a^2)=0,\] which implies that
\begin{equation}\label{Eq:a^2}
a^2=8dz(dz^2-1).
\end{equation} Since $dz$ is odd and $\gcd(dz,8(dz^2-1))=1$, it follows that $\abs{dz}$ and $\abs{(dz^2-1)/2}$ are integer squares. Thus, $d\mid z$ since $d$ is squarefree, and so $z=dn$ for some odd integer $n$. Then $w=2z=2dn$, and we see from \eqref{Eq:dw} that
\[4(1\pm 2^{m-2})=3dw^2=12d^3n^2.\] It follows that
 \begin{equation}\label{Eq:d^3n^2}
 d^3n^2=\left\{
 \begin{array}{cl}
   \dfrac{1-2^{m-2}}{3}<0 & \mbox{if $2\mid m$}\\[1em]
   \dfrac{1+2^{m-2}}{3}>0 & \mbox{if $2\nmid m$,}
 \end{array}\right.
 \end{equation}
 and consequently, 
 \[\abs{\dfrac{dz^2-1}{2}}=\abs{\dfrac{d^3n^2-1}{2}}=\left\{
 \begin{array}{cl}
   \dfrac{2^{m-3}+1}{3} & \mbox{if $2\mid m$}\\[1em]
   \dfrac{2^{m-3}-1}{3} & \mbox{if $2\nmid m$.}
 \end{array}\right.\]

 Recall that $\abs{(dz^2-1)/2}\in \Z$ is a square.
  However, if $2\mid m$ with $m\ge 6$, then
 \[\dfrac{2^{m-3}+1}{3}\equiv 3 \pmod{4}\] is not a square. Thus, $m=4$ in this case, so that $d=-1$ and $n=\pm 1$ from \eqref{Eq:d^3n^2}. From \eqref{Eq:a^2}, we have that
 \begin{equation}\label{Eq:Newa^2}
 a^2=8dz(dz^2-1)=8d^2n(d^3n^2-1),
 \end{equation}
 which implies that $n=-1$ and $a=\pm 4$. Using a computer algebra system and Theorem \ref{Thm:JKS}, it is easy to verify that both of the trinomials $x^4\pm 4x^3-1$ have Galois group isomorphic to $D_4$, but neither one is monogenic.
   Similarly, if $2\nmid m$ with $m\ge 7$, then
 \[\dfrac{2^{m-3}-1}{3}\equiv 5 \pmod{8}\] is not a square. Thus, $m=5$, so that $d^3n^2=3$ from \eqref{Eq:d^3n^2} in this case, which is clearly impossible.
  We conclude therefore that there are no monogenic trinomials $f(x)$ with $d\not \in \D$.

Suppose then that $d\in \D=\{1,2\}$. For each of these values of $d$, we use the fact that $r(t)=0$ and consider the integral points on the elliptic curve
\[E:\ Y^2=X^3-4d^3X,\]
where $X=dt$ and $Y=ad^2$. When $d=1$, Sage gives that the integral points $(X,Y)$ on $E$ are
\[(X,Y)\in \{(0,0),(\pm 2,0)\},\]
which yields no solutions for $(a,d)$ since $Y=ad^2=0$ for each point $(X,Y)$. When $d=2$, Sage gives that the integral points $(X,Y)$ on $E$ are
\[(X,Y)\in \{(4,\pm 8),(0,0),(8,\pm 16), (9,\pm 21), (1352,\pm 49712)\}.\] In this case, we get the viable solutions
\[(a,d)\in \{(\pm 2,2), (\pm 4,2), (\pm 12428,2)\}.\] Checking these possibilities reveals that
there are precisely two monogenic trinomials $f(x)$, namely
 \[x^4-2x^3+2 \quad \mbox{and} \quad x^4+2x^3+2,\]
 when $r(x)$ is reducible, and $\Gal(f)\simeq D_4$ for both of these trinomials, which completes the proof of item \eqref{Main1 I1}.
Furthermore, when $r(x)$ is irreducible, items \eqref{Main1 I2} and \eqref{Main1 I3} follow immediately from Theorem \ref{Thm:KW}.
\end{proof}
\section{The Proof of Theorem \ref{Thm:Main2}}\label{Section:Main2Proof}
\begin{proof} 

  Assume first that $(a,d)=(4k,27k^4+1)$, where $k\in \Z$ such that $k\ne 0$ and $d$ is squarefree.
We claim that $f(x)$ is irreducible over $\Q$.  
 Using calculus, it is easy to see that
 \begin{equation}\label{Eq:deff}
  f(x)=x^4+4kx^3+27k^4+1
  \end{equation} has an absolute minimum value of 1 at $x=-3k$, which implies that all zeros of $f(x)$ are non-real. Hence, if $f(x)$ is reducible over $\Q$, then $f(x)=g_1(x)g_2(x)$, where
 \[g_1(x)=x^2+Ax+B\in \Z[x]\quad \mbox{and} \quad g_2(x)=x^2+Cx+D\in \Z[x].\] Thus,
\begin{equation}\label{Eq:ffactored}
 f(x)=x^4+(A+C)x^3+(AC+B+D)x^2+(AD+BC)x+BD.
\end{equation}
Since $f(-3k)=1$, we deduce that $g_1(-3k)=g_2(-3k)$. We add this equation to the set of equations derived from equating coefficients in \eqref{Eq:deff} and \eqref{Eq:ffactored}, and we use Maple to solve this system. Maple gives two solutions, both of which have $k=0$, which contradicts our assumption that $k\ne 0$.  Hence, $f(x)$ is irreducible over $\Q$.

Next, we use Theorem \ref{Thm:JKS} to prove that $f(x)$ is monogenic. An easy calculation shows that 
\begin{equation*}\label{DelA4}
 \Delta(f)=\delta d^2=2^8(27k^4+1)^2.
 \end{equation*}  
  Consider first the prime divisor $q=2$ of $\Delta(f)$. 
  Note that $2\mid a$. If $2\nmid k$, then $4\mid d$, contradicting the fact that $d$ is squarefree. Hence, $2\mid k$ and $d\equiv 1 \pmod{4}$. Thus, $2\nmid d$, and it follows that condition \eqref{JKS:I2} of Theorem \ref{Thm:JKS} is satisfied since $2\mid a_2$ and $2\nmid d_1$.
  Next, suppose that $q\ne 2$ is a prime divisor of $\Delta(f)$ so that $q\mid d$.  Since $d$ is squarefree, we have that $q^2\nmid d$, so that condition \eqref{JKS:I1} of Theorem \ref{Thm:JKS} is satisfied. Thus, $f(x)$ is monogenic, and since $\delta=2^8$ is a square, we have that $\Gal(f)\simeq A_4$ by Theorem \ref{Thm:KW}, which establishes the theorem in this direction.

For the converse, assume that $f(x)$ is monogenic with $\Gal(f)\simeq A_4$.    Then $d$ is squarefree by Corollary \ref{Cor:squarefree}, and $\delta$ is a square by Theorem \ref{Thm:KW}.  Furthermore, if $d<0$, then $\delta<0$, which is impossible; and if $d=1$, then it is easy to check that $\delta$ is not a square if $a\ne 0$. Hence, $d\ge 2$. 
Suppose there is an odd prime $q$ such that $q\mid \delta$ and $q\nmid d$. Then $q\nmid a$. However, $q^2\mid \delta$ since $\delta$ is a square, which yields the contradiction that $f(x)$ is not monogenic by condition \eqref{JKS:I3} of Theorem \ref{Thm:JKS}.
Now suppose that there is an odd prime $q$ such that $q\mid \delta$ and $q\mid d$. Then $q$ divides $256d-\delta=27a^4$, so that $q^2\mid 27a^4$. Moreover, $q^2\mid \delta$ since $\delta$ is a square, and hence we have that $q^2$ divides $\delta+27a^4=256d$, contradicting the fact that $d$ is squarefree. Thus, we have shown that $\delta=2^{2m}$ for some integer $m\ge 0$. That is,
\begin{equation}\label{delta}
2^{2m}+27a^4=2^8d.
\end{equation} It is easy to verify that \eqref{delta} is impossible modulo 128 if $m\le 3$. Hence, $m\ge 4$ and $4\mid a$. Let $a=4k$ for some integer $k\ne 0$. We claim that $2\nmid d$. To establish this claim, we examine the exponent on the power of 2, denoted $\nu_2(*)$, that divides each side of the equation in \eqref{delta}. Since
\[\nu_2(a^4)=\nu_2(2^8k^4)=4z,\] for some integer $z\ge 2$, a straightforward computation shows that
 \begin{equation}\label{Eq:nu}
 \nu_2(2^{2m}+27a^4)=\left\{\begin{array}{cl}
   2m & \mbox{if $2m<4z$}\\
   4z & \mbox{if $2m>4z$}\\
   2^{2m+2}  & \mbox{if $2m=4z$}.\\
 \end{array}\right.
 \end{equation} If $2\mid d$, then $\nu_2(2^8d)=9$ since $d$ is squarefree, which contradicts \eqref{Eq:nu}. Hence, $2\nmid d$. With $q=2$, we see that if  $d\equiv 3 \pmod{4}$, then  condition \eqref{JKS:I2} of Theorem \ref{Thm:JKS} is not satisfied since $2\mid a_2$ and $2\mid d_1$, which contradicts the fact that $f(x)$ is monogenic. Thus, $d\equiv 1 \pmod{4}$.
Suppose now that $m\ge 5$. Hence, it follows from \eqref{delta} that
\[27k^4=d-2^{2m-8}\equiv 1 \pmod{4},\] which is impossible, since $27k^4 \Mod{4}\in \{0,3\}$. Therefore, $m=4$ and it is easy to see from \eqref{delta}, with $a=4k$, that $d=27k^4+1$, which completes the proof of the theorem.
\end{proof}

\section*{Acknowledgments} The authors thank the anonymous referee, especially for catching some careless mistakes in the proof of Theorem \ref{Thm:Main1}.





\end{document}